\documentclass{amsart}
\usepackage{mathrsfs}
\usepackage{amsfonts,amssymb,amsmath}
\def\fk{\mathfrak}
\def\NN{{\mathbb N}}
\def\ZZ{{\mathbb Z}}
\def\RR{{\mathbb R}}

\def\eps{\varepsilon}
\def\al{\alpha}

\def\lm{\lambda}
\def\Lm{\Lambda}

\def\rho{\varrho}
\def\phi{\varphi}
\def\gm{\gamma}

\def\Ups{\Upsilon}

\newcommand{\f}[2]{\frac{#1}{#2}}
\newcommand{\be}[1]{\begin{equation}\label{#1}}
\newcommand{\ee}{\end{equation}}
\newcommand{\multsum}[2]{\sum_{{\scriptstyle #1}\atop {\scriptstyle #2}}}
\newcommand{\trisum}[3]{\sum_{{{\scriptstyle #1}\atop {\scriptstyle #2}}\atop{\scriptstyle #3}}}

\newcommand{\rf}[1]{{\rm (\ref{#1})}}

\def\cal{\mathscr}

\def\dlt{\delta}

\def\Ups{\Upsilon}

\def\fk{\mathfrak}
\def\Ups{\Upsilon}

\def\bfx{{\bf x}}\def\bfy{{\bf y}}

\newcommand{\bflm}{\mbox{\boldmath$\lambda$}}
\newcommand{\bflms}{\mbox{\boldmath$\scriptstyle\lambda$}}
\newcommand{\bfLm}{\mbox{\boldmath$\Lambda$}}
\newcommand{\bfsg}{\mbox{\boldmath$\sigma$}}
\newcommand{\bfsgs}{\mbox{\boldmath$\scriptstyle\sigma$}}
\newcommand{\bfbe}{\mbox{\boldmath$\beta$}}

\newtheorem*{thm}{Theorem}
\newtheorem{lem}{Lemma}

\begin{document}

\title{Random diophantine inequalities of additive type}
\author{J\"org Br\"udern and Rainer Dietmann}
\address{
J\"org Br\"udern, Mathematisches Institut, Bunsenstrasse 3--5, 37073 G\"ottingen, Germany.}
\email{bruedern@uni-math.gwdg.de}
\address{
Rainer Dietmann, Department of Mathematics, Royal Holloway, University of London,
Egham, Surrey TW20 0EX, UK.}
\email{Rainer.Dietmann@rhul.ac.uk}
\thanks{Rainer Dietmann acknowledges support by EPSRC grant EP/I018824/1 'Forms in many variables'}
\begin{abstract}
Using the Davenport-Heilbronn circle method, we show that for almost all
additive Diophantine inequalities of degree $k$ in more than $2k$ variables
the expected asymptotic formula for the density of solutions holds true.
This appears to be the first metric result on Diophantine inequalities.
\end{abstract}
\maketitle

\section{Introduction}
The study of additive diophantine inequalities has been one of the guiding 
themes in diophantine approximation. Let $k$ and $s$ be natural 
numbers, $s\ge 2$. For non-zero real numbers $\lambda_1,\ldots,\lambda_s$ 
and a positive number $\tau$, consider the inequality
\begin{equation}
|\lambda_1x^k_1+\lambda_2x^k_2+\ldots+\lambda_sx^k_s|<\tau\label{1}
\end{equation}
that is to be solved in integers $x_j$, not all zero. 
Leaving aside the most classical linear case $k=1$ with an overwhelmingly 
rich literature (see \cite{HL}, \cite{sbook}, \cite{B} and the 
references therein for themes related to this article), already the case $k=2$ shows typical features. There are 
two obstacles to solubility. On the one hand, the quadratic form on the left 
hand side of \rf{1} may be definite in which case the only solution of \rf{1} 
is $x_1=\ldots=x_s=0$, at least when $\tau$ is small. On the other hand, the 
form may be indefinite, but a real multiple of a form with integer coefficients.
 Then, again for small $\tau$, the inequality \rf{1} is satisfied if and only 
if the form with integer coefficients vanishes, and when $s=3$ or $4$, there 
may be $p$-adic obstructions to realize this. The remaining cases are described 
by the conditions that $\lambda_1,\ldots,\lm_s$ are not all of the same sign, and 
that at least one of the ratios $\lambda_i/\lambda_j$ is irrational, and when 
$s\ge 3$, one would expect non-trivial solutions of \rf{1} to exist. 
This surpasses a long-standing conjecture of Oppenheim, but nowadays is merely a special case of celebrated work of Margulis \cite{M}.

Similar results are expected for larger values of $k$. The potential obstacles 
to solubility are the same, but ``definite'' forms exist only for even degree. 
Therefore, we write $\Lm^{(s)} = (\RR\setminus\{0\})^s$ and
 $\Lambda^{(s)}_k=\Lambda^{(s)}$ when $k$ is odd, but put
$$
\Lambda^{(s)}_k=\{\bflm\in\Lambda^{(s)}:\mbox{ $\lambda_i\lambda_j<0$ for 
some $1\le i,j\le s$}\}
$$
when $k$ is even. 

Davenport and Heilbronn \cite{DH} showed that when $s>2^k$, $ 
\bflm\in\Lambda^{(s)}_k$ and some ratio $\lambda_i/\lambda_j$ is irrational, 
then \rf{1} admits infinitely many solutions ${\bf x}\in\mathbb{Z}^s$. Their 
pivotal contribution was very influential, and the Fourier transform method 
that they developed still underpins much recent work. A first wave of 
refinements led to a reduction of the variables required, and it was also 
realized that one could take $\tau=|{\bf x}|^{-\sigma}$, with 
$|{\bf x}|=\max|x_j|$ and some suitably small $\sigma>0$, and still 
guarantee the 
existence of infinitely many integer solutions of \rf{1}; see 
\cite{hlm}, chapter 11.

A major innovation is due to Freeman \cite{F1, F2}. Inspired by related 
work of Bentkus and G\"otze \cite{BG}, he considered the number 
$N(P)=N^{(k)}_{\bflms}(P;\tau)$ of integer solutions of \rf{1} within the box 
$|{\bf x}|\le P$. When $s>2^k$  or $s\ge (1+o(1))k^2 \log k $, and 
$\bflm\in\Lambda^{(s)}_k $ has some  $\lambda_i/\lambda_j$ irrational, Freeman
 established that the limit
\begin{equation}
\lim\limits_{P\to\infty}P^{k-s}N(P)\label{2}
\end{equation}
exists and is positive. Wooley \cite{W1} has smaller admissible values 
for $s$, and on combining Theorem 1.1 of \cite{W1}  with Wooley's very recent 
furious work \cite{W2,W3}, one obtains Freeman's result under the
less restrictive condition $s\ge 2k^2$.

Perhaps Freeman's asymptotic formula \rf{2} remains valid for $s>k$, 
but the range $k<s\le 2k$ has resisted attacks even subject to the strongest 
plausible hypotheses on Weyl sums. This is due to the familiar square-root 
cancellation barrier for estimates of exponential sums. When $s>2k$, however, 
we are able to demonstrate an asymptotic formula for $N^{(k)}_{\bflms}(P;\tau)$, 
with a strong error term and uniform with respect to $\tau$, for almost all 
$\bflm\in\Lm^{(s)}_k$ (in the sense of Lebesgue measure).
\begin{thm}
Let $k\ge2$, $s>2k$ and $\delta=8^{-2k}$. Then, for almost all ${\bfLm}^{(s)}_k$, 
there exist a number $P_0=P_0(\bflm,k)$, and a positive real number 
$J_{k,s}(\bflm)$ such that the inequality
\begin{equation}
|N^{(k)}_{\bflms}(P;\tau)-2\tau J_{k,s}(\bflm)P^{s-k}|<P^{s-k-\delta}\label{3}
\end{equation}
holds for all $P\ge P_0$, uniformly in $0<\tau\le 1$.
\end{thm}
No effort has been made to optimize the value of $\delta$. It is expected 
that it is possible to construct $\bflm\in\bfLm^{(s)}_k$ 
where the convergence in \rf{2} is slower than any predetermined speed.  
Therefore, it is  rather remarkable that one saves a fixed power of $P$ in \rf{3}, 
outside a set of Lebesgue measure $0$. This much is new even in the 
case $k=2$ where the condition $s>2k$ coincides with Freeman's $s>2^k$.

Our approach depends on some simple observations concerning differences of 
two integral $k$-th powers that we collect in the next section. The results are 
then applied within the estimation of a mean value for Weyl sums. This is 
the theme of section 3. In some sense, Lemma 4 below may 
be considered as an averaged version of Hardy and Littlewood's famous 
conjecture K. Indeed, the conclusion of the theorem would be valid for {\em all}
$\bflm\in\Lm_k^{(s)}$ if the truth of hypothesis K were postulated. With the main lemma now in hand, we
proceed to set the scene for the Fourier transform 
method and explain the strategy of proof in the short section 4. The singular 
integral $J_{k,s}(\bflm)$ will be discussed in detail in section 5. Then, the 
classical interference  principle will be studied from a new perspective,
providing a suitable estimate for the complementary compositum 
in our application of the Davenport-Heilbronn method.
This is the theme of section 6. In section 7, we obtain a weighted version of 
our theorem, and in the final section, the weights will be removed to 
complete the proof. 

\section{Differences of two k-th powers} For the rest of this article,
suppose that $k$
is a  natural number, $k\ge 2$.
Let $r\in\mathbb N$, 
and let $S_r(P,Z)$ be the number of $\bfx,\bfy\in\mathbb Z^r$ with 
$|\bfx|\le P$, $|\bfy|\le P$ and 
$$
1\le|x^k_r-y^k_r|\le|x^k_{r-1}-y^k_{r-1}|\le\ldots\le|x^k_1-y^k_1|\le Z.
$$
\begin{lem}
 One has
$$S_r (P,Z)\ll Z^{r/(k-1)}P^{r(k-2)/(k-1)}.$$
\end{lem}
\noindent
{\it Proof.} First consider the case where $r=1$ and $k$ is even. By symmetry, 
it suffices to bound the number of pairs $x,y$ with $1\le x^k-y^k\le Z$ and 
$0\le y<x\le P$. Then, the integer $h=x-y$ satisfies $1\le h\le P$ and 
$Z\ge x^k-y^k\ge hx^{k-1}$ so that
$$
S_1(P,Z)\ll\sum\limits_{1\le h\le P}\Big(\frac{Z}{h}\Big)^{1/(k-1)}\ll 
Z^{1/(k-1)}P^{(k-2)/(k-1)},
$$
as required.

Next, suppose that $r=1$, and that $k$ is odd. The number of pairs $x_1,y_1$ 
with $x_1y_1\ge 0$ that are counted by $S_1(P,Z)$ 
can be estimated exactly as in the case where $k$ was even. For the remaining 
pairs, one observes symmetry, and it will suffice to count those where 
$x_1>0$, $y_1<0$. Then $x^k_1+|y_1^k|\le Z$, and a rough lattice point count 
shows that there are no more than $O(Z^{2/ k})$ such pairs of integers. This 
suffices to confirm the claim of the lemma when $Z\le P^k$, and when 
$Z\ge P^k$, the estimate proposed in Lemma 1 is weaker than the obvious bound 
$S_1(P,Z)\ll P^2$. This establishes Lemma 1 in the case $r=1$.

Now let $r\ge 1$, and note that
$$
S_{r+1}(P,Z)=
\multsum{|x_1|\le P,|y_1|\le P}{1\le |x^k_1-y^k_1|\le Z}S_{r}(P,|x^k_1-y^k_1|).
$$
We proceed by induction on $r$ and estimate $S_r$ by the 
induction hypothesis. For any 
$x_1,y_1$ occurring in the sum, we may then suppose that
$$
S_{r}(P,|x^k_1-y^k_1|)\ll Z^{r/(k-1)}P^{r(k-2)/(k-1)}.
$$
Hence,
$$
S_{r+1}(P,Z)\ll Z^{r/(k-1)}P^{r(k-2)/(k-1)}S_1(P,Z).
$$
Using the bound for $S_1(P,Z)$ that we have already established, one completes 
the induction and the proof of Lemma 1.

\begin{lem} Let $\beta$ be a real number with $1/(k-1)\le \beta\le 1$. Then
$$
\multsum{|x|\le P,|y|\le P}{x^k\neq y^k}|x^k-y^k|^{-\beta}\ll P^{1-\beta}\log P.
$$
\end{lem}
The proof is very similar to the proof of Lemma 1. When $k$ is even, it again 
suffices to estimate the portion of the sum where $0\le y<x\le P$. For $h=x-y$ 
one finds $P^k\ge x^k-y^k\ge hx^{k-1}$, so that
$$
\multsum{|x|\le P,|y|\le P}{x^k\neq y^k}|x^k-y^k|^{-\beta}
\ll\sum\limits_{1\le h\le P} h^{-\beta}\sum_{1\le x\le P}
x^{-\beta(k-1)}\ll P^{1-\beta}\log P.
$$

When $k$ is odd, the contribution of pairs $x$, $y$ with $xy\ge 0$  can be estimated as above. It remains to consider the portion of the sum 
in question where $x$ and $y$ have opposite signs. In that case, one has 
$|x^k-y^k|\ge \max(|x|^k,|y|^k)$, so that this portion of the sum does not 
exceed
$$
\ll\sum\limits_{1\le|x|\le P}\sum\limits_{|y|\le|x|}|x|^{-\beta k}
\ll\sum\limits_{1\le x\le P}x^{1-\beta k}.
$$
This expression is bounded by $O(1)$ when $\beta>2/ k$, and is 
$O(P^{2-\beta k}\log P)$ for $1/(k-1)\le\beta\le 2/ k$. These bounds are 
stronger than required to complete the proof of Lemma 2. 

For $r\in\mathbb N$, let $\cal S_r(P)$ denote the set of $\bfx,\bfy\in\mathbb Z^r$ with $|\bfx|\le P,|\bfy|\le P$ and
$$
|x^k_1-y^k_1|\ge |x^k_2-y^k_2|\ge\ldots\ge |x^k_r-y^k_r|\ge 1.
$$
\begin{lem}
 Let $1\le r\le k$. Then
$$
\sum\limits_{(\bfx,\bfy)\in\cal S_r(P)}|x^k_1-y^k_1|^{-1}\ll P^r\log P.
$$
\end{lem}
\noindent
{\it Proof.} By Lemma 1, the sum in question does not exceed
$$
\trisum{|x_1|\le P}{|y_1|\le P}{x^k_1\neq y^k_1}  
\f{S_{r-1}(P,|x^k_1-y^k_1|)}{ |x^k_1-y^k_1|} \ll P^{(r-1)(k-2)/(k-1)}
\trisum{|x_1|\le P}{|y_1|\le P}{x^k_1\neq y^k_1}|x^k_1-y^k_1|^{(r-1)/(k-1)-1}.
$$
When $1\le r\le k-1$, Lemma 2 is applicable to the sum on the right hand side, and provides the desired estimate. When $r=k$, the sum over $x_1,y_1$ on the right is $O(P^2)$, and the estimate proposed in Lemma 3 again follows.

\section{The catalytic mean value} The main auxiliary estimate concerns a certain mean value. Let $C$ be a fixed real number with $C\ge 1$. For $0<\eta\le 1$, define a measure $d_\eta\alpha$ on $\mathbb R$ by
\begin{equation}
d_\eta\alpha=\eta\Big(\frac{\sin\pi\eta\alpha}{\pi\eta\alpha}\Big)^2\,d\alpha
\label{4}
\end{equation}
where $d\alpha$ is the standard Lebesgue measure. Its Fourier transform is
\begin{equation}
W_\eta(\alpha)=\int^\infty_{-\infty}e(-\alpha\beta)\,d_\eta\beta
=\max\Big(0,1-\frac{|\alpha|}{\eta}\Big).\label{5}
\end{equation}
We introduce the Weyl sum
\begin{equation}
f(\alpha)=\sum\limits_{|x|\le P}e(\alpha x^k)\label{6}
\end{equation}
and consider the moment
$$
\Xi=\int_{|\bflms|\le C}\int^\infty_{-\infty}|f(\lambda_1\alpha)f(\lambda_2\alpha)\ldots f(\lm_{2k}\alpha)|\,d_\eta\alpha \,d\bflm.
$$
\begin{lem}
 Let $C\ge 1$. Then, uniformly for $0<\eta\le 2$, one has
$$
\Xi\ll P^k\log P.
$$
\end{lem}
\noindent
{\it Proof.} By the Fubini-Tonelli theorem and Schwarz's inequality,
\begin{eqnarray*}
\Xi&=&\int^\infty_{-\infty}\Big(\int^C_{-C}|f(\lm\alpha)|\,d\lm\Big)^{2k}\,d_\eta\alpha\\
&\le&(2C)^k\int^\infty_{-\infty}
\Big(\int^C_{-C}|f(\lm\alpha)|^2\,d\lm\Big)^k\,d_\eta\alpha.
\end{eqnarray*}
Now reverse the order of integration again. By \rf{5} and \rf{6}, this yields
$$
\Xi\le(2C)^k\multsum{|\bfx|\le P}{|\bfy|\le P}\int_{[-C,C]^k}W_\eta(\lm_1(x^k_1-y^k_1)+\ldots+\lm_k(x^k_k-y^k_k))\,d\bflm.
$$
By symmetry, it suffices to estimate the portion of the sum on the right hand side where
$$
|x^k_1-y^k_1|\ge|x^k_2-y^k_2|\ge\ldots\ge|x^k_k-y^k_k|.
$$

Subject to this additional constraint, first consider the contribution of terms with $x^k_1=y^k_1$. Then $x^k_j=y^k_j$ for all $1\le j\le k$, leaving $O(P^k)$ choices for $\bfx,\bfy$. For any such choice, the integrand is $W_\eta(0)=1$. Consequently, the contribution of these terms to $\Xi$ is $O(P^k)$.

For the remaining terms, there is a number $r$ with $1\le r\le k$ and such that $|x^k_r-y^k_r|\ge 1$ but $x^k_j=y^k_j$ for $j>r$ (if any). 
Then there are no more than $(2P+1)^{k-r}$ 
choices for $x_j,y_j$ with $r<j\le k$, and the integrand in the penultimate display does not depend on $\lm_{r+1},\ldots,\lm_k$. Hence, the contribution to $\Xi$ that arises from terms $\bfx,\bfy$ with a fixed value of $r$ does not exceed
$$
(2C)^{2k-r}(2P+1)^{k-r}\sum\limits_{(\bfx,\bfy)\in\cal S_r(P)}\int_{[-C,C]^r}W_\eta(\lm_1(x^k_1-y^k_1)+\ldots+\lm_r(x^k_r-y^k_r))\,d\bflm.
$$
One integrates over $\lm_1$ first, considering $(\bfx,\bfy)\in\cal S_r(P)$ and $\lm_2,\ldots,\lm_r\in\mathbb R$ as fixed real numbers. By \rf{5}, the integrand is non-zero only on an interval for $\lm_1$, of length $2\tau/|x^k_1-y^k_1|$, and one has $0\le W_\eta\le 1$. It follows that the expression in the previous display does not exceed
$$
\ll P^{k-r}\sum\limits_{(\bfx,\bfy)\in\cal S_r(P)}\frac{2\tau}{|x^k_1-y^k_1|}\ll P^k\log P.
$$
Here, Lemma 3 was applied to confirm the rightmost inequality. The lemma follows by summing the various contributions.

\section{The Fourier transform method} We prepare the scene for an application of the Davenport-Heilbronn method, as renovated in \cite{ARTS8}. Let $0<\eta\le 2$, and let $\bflm\in\Lm^{(s)}_k$. Consider the weighted analogue of $N^{(k)}_{\bflms}(P,\eta)$ defined by
\begin{equation}
I_{\bflms}(P,\eta)=\sum\limits_{|\bfx|\le P}W_\eta(\lm_1x^k_1+\lm_2x^k_2+\ldots+\lm_sx^k_s)\label{7}
\end{equation}
where $W_\eta$ is the function defined in \rf{5}. By \rf{6}, we have the alternative representation
\begin{equation}
I_{\bflms}(P,\eta)=\int^\infty_{-\infty}f(\lm_1\alpha)f(\lm_2\alpha)\ldots f(\lm_s\alpha)\,d_\eta\alpha.\label{8}
\end{equation}

We shall derive an asymptotic formula for this integral 
that will hold for almost all $\bflm\in\Lm^{(s)}_k$. The main term in this formula arises from the central interval
\begin{equation}
\mathfrak{C}=[-P^{1/3-k},P^{1/3-k}],\label{9}
\end{equation}
and the contribution from the complementary compositum
\begin{equation}
\mathfrak{c}=\{\alpha\in\mathbb{R}:|\alpha|>P^{1/3-k}\}\label{10}
\end{equation}
will be negligible on average over $\bflm$.

From now on, suppose that $k,s$ with $k\ge2$, $s>2k$, are fixed once and for all. Also, let $C$ be a real number with $C\ge2$, and let
\begin{equation}
\Lm^{(s)}_k(C)=\{\bflm\in\Lm^{(s)}_k:1/C\le\lm_j\le C\;(1\le j\le s)\}.\label{11}
\end{equation}
Implicit constants in estimates to follow will depend on the parameters $s,k,C$, but are  uniform with respect to $P,\eta$ and $\bflm\in\Lm^{(s)}_k(C)$.

When $\mathfrak{A}\subset\mathbb R$ is a measurable set, let
\begin{equation}
\cal I(\mathfrak A)=\int_{\mathfrak A}f(\lm_1\alpha)f(\lm_2\alpha)\ldots f(\lm_s\alpha)\,d_\eta\alpha.\label{12}
\end{equation}
In the interest of notational compactness, dependence of $\cal{I}(\mathfrak{A})$ on $\bflm,\eta,P$ has been suppressed. Note that by \rf{8}, \rf{9}, \rf{10} and \rf{12},
one has
\be{12a} I_{\bflms}(P,\eta)= \cal{I}(\mathbb{R}) = \cal I(\mathfrak C)+\cal I(\mathfrak c).\ee

\section{The central interval} An asymptotic formula will be provided for $\cal{I}(\mathfrak{C})$. The argument is largely standard, but there is no appropriate reference for the uniformity issue relevant for the current considerations. We therefore indulge into a detailed account, but we shall be brief. The exposition is modelled on Wooley \cite{W1} where appropriate, but there are differences because precise control of error terms is needed.

Let
\begin{equation}
v(\alpha)=\int^P_{-P}e(\alpha\xi^k)\,d\xi.\label{13}
\end{equation}
By \rf{6} and partial summation,
$$
f(\alpha)=v(\alpha)+O(1+P^k|\alpha|)
$$
uniformly for $\alpha\in\mathbb{R}$, and partial integration readily yields
$$
v(\alpha)\ll P(1+P^k|\alpha|)^{-1/k}.
$$
Now let $\bflm\in\Lm^{(s)}_k(C)$, and put
$$
V_{\bflms}(\alpha)=v(\lm_1\alpha)v(\lm_2\alpha)\ldots v(\lm_s\alpha).
$$
Then, by the preceding estimates, for $\alpha\in\mathfrak{C}$, one has
$$
f(\lm_1\alpha)f(\lm_2\alpha)\ldots f(\lm_s\alpha)=V_{\bflms}(\alpha)+O(P^{s-2/3}).
$$
We integrate over $\mathfrak{C}$, against $d_\eta\alpha$. By \rf{4}, \rf{9} and \rf{12}, this implies that uniformly in $0<\eta\le 2$, $\bflm\in\Lm^{(s)}_k(C)$, one has
$$
\cal{I}(\mathfrak{C})=\int_{\mathfrak{C}}V_{\bflms}(\alpha)\,d_\eta\alpha+O(P^{s-k-1/3}).
$$
Now use $s>2k$ and
the upper bound for $v(\alpha)$ to infer that 
for $\bflm\in\Lm^{(s)}_k(C)$, one has $V_{\bflms}(\alpha)\ll P^s(1+P^k|\alpha|)^{-2}$. Whence
$$
\int_{\mathfrak{c}}|V_{\bflms}(\alpha)|\,d_\eta\alpha\ll P^{s-k-1/3},
$$
and one may add this to the previous display to conclude that
\begin{equation}
\cal{I}(\mathfrak{C})=\int^{\infty}_{-\infty}V_{\bflms}(\alpha)\,d_\eta\alpha+O(P^{s-k-1/3})\label{14}
\end{equation}
holds uniformly for $0<\eta\le 2$, $\bflm\in\Lm^{(s)}_k(C)$.

It remains to evaluate the integral on the right hand side of \rf{14}. It will be helpful to write $V_{\bflms}(\alpha)$ as a Fourier transform. Let $\bfsg=(\sigma_1,\ldots,\sigma_s)$ with $\sigma_j=\pm 1$ for $1\le j\le s$, and let $\bflm\in\Lm^{(s)}_k(C)$. Then, for a real parameter $\beta$, let
$
\cal{B}(\beta,\bflm)$
be the set of all $(\beta_2,\ldots,\beta_s)\in\RR^{s-1}$ satisfying the inequalities
$$
0\le\beta_j\le|\lm_j|\quad(2\le j\le s),\quad0\le\beta-\sigma_1\sigma_2\beta_2-\ldots-\sigma_1\sigma_s\beta_s
\le|\lm_1|.
$$
Then, the integral
\begin{equation}
E_{\bfsgs}(\beta)= \int_{\cal{B}(\beta,\bflms)}(\beta-\sigma_1\sigma_2\beta_2-\ldots- \sigma_1\sigma_s\beta_s)^{\frac{1}{k}-1} (\beta_2\beta_3\ldots\beta_s)^{\frac{1}{k}-1}\,d(\beta_2,\ldots,\beta_s)\label{15}
\end{equation}
defines  a non-negative, continuous and compactly supported function.

Beyond this point, the details of the argument depend on the parity of $k$. Hence, we temporarily restrict our attention to the case where $k$ is even. We take $\sigma_j=\lm_j/|\lm_j|$ and note that for $\bflm\in \Lm^{(s)}_k(C)$, not all $\sigma_j$ are equal. An obvious substitution yields
$$
v(\lm_j\alpha)=2\int^P_0e(\lm_j\alpha\xi^k)\,d\xi
=\frac{2P}{k}|\lm_j|^{-1/k}\int^{|\lm_j|}_0 \beta^{1/k-1}e(\alpha\beta\sigma_jP^k)\,d\beta,
$$
and so,
$$
V_{\bflms}(\alpha)=
\Big(\frac{2P}{k}\Big)^s|\lm_1\ldots\lm_s|^{-1/k}
\int_{\cal{U}(\bflms)}(\beta_1\ldots\beta_s)^{1/k-1}e(\alpha P^k(\sigma_1\beta_1+\ldots+\sigma_s\beta_s))\,d\bfbe
$$
where $\cal{U}(\bflm)=[0,|\lm_1|]\times \ldots \times[0,|\lm_s|]$. We substitute $\beta$ for $\beta_1$ via
$$
\sigma_1\beta=\sigma_1\beta_1+\sigma_2\beta_2+\ldots+\sigma_s\beta_s
$$
to see that
\begin{equation}
V_{\bflms}(\alpha)=\Big( \f{2P}{k}\Big)^s|\lm_1\ldots\lm_s|^{-1/k}\int^\infty_{-\infty}E_{\bfsgs}(\beta)e(\sigma_1\alpha P^k\beta)\,d\beta.\label{16}
\end{equation}
By Fourier's integral theorem, we then have
\begin{equation}
E_{\bfsgs}(\beta)=\Big(\frac{k}{2P}\Big)^s|\lm_1\ldots\lm_s|^{1/k}\int^\infty_{-\infty}V_{\bflms}(\alpha P^{-k})e(-\sigma_1\alpha\beta)\,d\alpha.\label{17}
\end{equation}
By \rf{16}, Fubini's theorem and \rf{5},
$$
\int^\infty_{-\infty}V_{\bflms}(\alpha)\,d_\eta\alpha=
\Big(\frac{2P}{k}\Big)^s|\lm_1\ldots\lm_s|^{-1/k}\int^\infty_{-\infty}E_{\bfsgs}(\beta)W_\eta(\beta P^k)\,d\beta.
$$
Note that the integrand on the right hand side vanishes unless $|\beta|\le P^{-k}$. In the latter range for $\beta$, one may use \rf{17} to conclude that for $\bflm\in\Lm^{(s)}_k(C)$ one has
\begin{eqnarray*}
E_{\bfsgs}(\beta)-E_{\bfsgs}(0)&=&
\int^\infty_{-\infty}\Big(\frac{k}{2P}\Big)^s|\lm_1\ldots\lm_s|^{1/k}V_{\bflms}(\alpha P^{-k})(e(-\sigma_1\alpha\beta)-1)\,d\alpha\\
&\ll&\int^\infty_{-\infty}(1+|\alpha|)^{-s/k}|e(\alpha\beta)-1|\,d\alpha\\
&\ll&P^{-k},
\end{eqnarray*}
since $s>2k$.
Hence, in the preceeding identity, $E_{\bfsgs}(\beta)$ may be replaced by $E_{\bfsgs}(0)$, and one concludes that
$$
\int^\infty_{-\infty}V_{\bflms}(\alpha)\,d_\eta\alpha=\eta J_{k,s}(\bflm)P^{s-k}+O(P^{s-2k})
$$
where
$$
J_{k,s}(\bflm)=\Big(\frac{2}{k}\Big)^s|\lm_1\ldots\lm_s|^{-1/k}E_{\bfsgs}(0).
$$
By the definition of $E_{\bfsgs}(\beta)$, one finds that $E_{\bfsgs}(0)>0$ because the $\sigma_j$ are not all of the same sign, and one also sees that
$\Lm^{(s)}_k\to \RR$,      
$\bflm\mapsto E_{\bfsgs}(0)$ is continuous. By \rf{14}, we now deduce the case when $k$ is even of the following lemma.
\begin{lem}
 Let $s>2k$, $k\ge2$ and $C\ge2$. Then, there exists a continuous function $J_{k,s}:\Lm^{(s)}_k\to(0,\infty)$ with the property that
$$
\cal{I}(\mathfrak{C})=\eta J_{k,s}(\bflm)P^{s-k}+O(P^{s-k-1/3})
$$
holds uniformly for $0<\eta\le 2$ and $\bflm\in\Lm^{(s)}_k(C)$.
\end{lem}
It remains to establish this lemma when $k$ is odd. In this case,
$$
v(\lm_j\alpha,P)=\int^P_0(e(\lm_j\alpha\xi^k)+e(-\lm_j\alpha\xi^k))\,d\xi.
$$
Proceeding as before, the consequential analogue of \rf{16} is the identity
$$
V_{\bflms}(\alpha)=\Big(\frac{P}{k}\Big)^s|\lm_1\ldots\lm_s|^{-1/k}
\multsum{\sigma_j=\pm1}{j=1,\ldots,s}
\int^\infty_{-\infty}E_{\bfsgs}(\beta)e(\sigma_1\alpha P^k\beta)\,d\beta
$$
where we now consider $E_{\bfsgs}(\beta)$ as a function of the {\it independent} parameters $\sigma_j\in\{1,-1\}$ and $\bflm\in\Lm^{(s)}_k(C)$. One may then follow through the argument used in the case when $k$ is even. Lemma 5  follows with
$$
J_{k,s}(\bflm)=k^{-s}|\lm_1\ldots\lm_s|^{-1/k}\sum\limits_{\sigma_j=\pm1}E_{\bfsgs}(0),
$$
in which one has $E_{\bfsgs}(0)\ge0$, and $E_{\bfsgs}(0)>0$ whenever the $\sigma_j$ are not all of the same sign. 

\section{The interference estimate} 
In traditional applications of the Davenport-Heilbronn method, the
treatment of the complementary compositum depends on an interference
principle. This asserts that whenever $|\al|$ is neither too large 
nor too small, and $\lm_1/\lm_2$ is irrational, then the product
$|f(\lm_1\al)f(\lm_2\al)|$ is rather smaller than the trivial
bound $P^2$. We need this in strong quantitative form, but may take
advantage of averages over $\lm_j$. 

\begin{lem} Let $k\ge 2$, $s>2k$ and $C\ge 2$. Then, uniformly in
$0<\eta\le 2$, one has
$$ \int_{\Lm_k^{(s)}} \int_{\fk c} |f(\lm_1\al)\ldots f(\lm_s\al)|\,\,d_\eta\al
\,d\bflm \ll P^{s-k-10\dlt}. $$
\end{lem}
\noindent
{\em Proof}. Consider the tail ${\fk t}=\{\al: |\al|\ge 1\}$, and when
$1\le j\le s$ and $\bflm\in\Lm_k^{(s)}$, let
$$ {\fk d}_j(\bflm) = \{P^{1/3-k}\le |\al|\le 1: |f(\lm_j\al)|\le P^{1-11\dlt}\}.$$
Furthermore, let
$$ {\fk D}(\bflm) = \{P^{1/3-k}\le |\al|\le 1: |f(\lm_j\al)|> P^{1-11\dlt}
\, (1\le j\le s)\}. $$
Then $\fk c$ is the union of the sets $\fk t$, ${\fk D}(\bflm)$ and the
$s$ sets ${\fk d}_j(\bflm)$. Thus, on writing
$$ F(\al) = |f(\lm_1\al)\dots f(\lm_s\al)| $$
in the interest of brevity, the integrals 
$$ \Ups^* = \int_{\Lm_k^{(s)}(C)} \int_{\fk t} F(\al) \,d_\eta\al\,d\bflm,$$
$$ \Ups_0 =  \int_{\Lm_k^{(s)}(C)} \int_{{\fk D}(\bflms)} F(\al) \,d_\eta\al\,d\bflm,$$
and, for $1\le j\le s$, 
$$ \Ups_j =  \int_{\Lm_k^{(s)}(C)} 
\int_{{\fk d}_j(\bflms)} F(\al) \,d_\eta\al\,d\bflm$$
all exist by Tonelli's theorem, and one has
\be{A} \int_{\Lm_k^{(s)}} \int_{\fk c} F(\al) \,d_\eta\al\,d\bflm
\le \Ups^* + \sum_{j=0}^s \Ups_j.\ee

The estimation of $\Ups^*$ is straightforward. First observe that, by orthogonality,
\be{H1} \int_0^1 |f(\gm)|^2 \,d\gm \ll P. \ee
The main argument begins much as in the
proof of Lemma 4 with an exchange of the order of integration. Then, by 
Fubini's theorem and symmetry, 
$$  \Ups^* \le 2 \int_1^\infty \Big(\int_{-C}^C |f(\lm\al)|\,d\lm\Big)^s 
\,d_\eta\al. $$
Now suppose that $\al\ge 1$. By Schwarz's inequality and an obvious
substitution,
$$\int_{-C}^C |f(\lm\al)|\,d\lm \le (C/\al)^{1/2} \Big(\int_{-\al C}^{\al C}
|f(\gm)|^2\,d\gm\Big)^{1/2}. $$
But $\al C\ge 2$, and $f(\gm)$ has period 1. Hence, by \rf{H1}, the bound
$$ \int_{-C}^C |f(\lm\al)|\, d\lm \ll P^{1/2} $$
holds uniformly for $\al\ge 1$, with an implicit constant depending only on $C$.
Consequently, uniformly in $0<\eta \le 2$, one finds that
$$ \Ups^* \ll P^{s/2} \int_{-\infty}^\infty \,d_\eta\al \ll P^{s/2}. $$

Now consider $\Ups_1$. By the definition of ${\fk d}_1(\bflm)$,
$$ \Ups_1 \le P^{1-11\dlt} \int_{|\bflms|\le C} \int_{-\infty}^\infty
|f(\lm_2\al)\ldots f(\lm_s\al)|\,d_\eta\al \,d\bflm. $$
The condition that $s>2k$ assures that at least $2k$ of the $\lm_j$
occur in the integrand on the right hand side. Hence, one may apply 
Lemma 4 and perform remaining integrations (if any) trivially. This yields
$$ \Ups_1 \ll P^{s-k-11\dlt}\log P. $$
By symmetry, the same bound is valid for $\Ups_j$ when $1\le j\le s$.

It remains to estimate $\Ups_0$.   
Reversing the order of integrations yields
$$ \Ups_0 = 2 \int_{P^{1/3-k}}^1 \int_{M_\al(C)} F(\al) \,d\bflm\,d_\eta\al $$
where
$$ M_\al(C) =\{ \bflm\in\Lm_k^{(s)}(C): \,|f(\lm_j\al)|\ge P^{1-11\dlt}\, 
(1\le j\le s)\}.
$$
 
We proceed to show that
uniformly for $\al$ in the range  $P^{1/3-k} <\al\le 1$, the measure of $M_\al(C)$ does not
exceed $O(P^{(55\dlt k -k)(s-1)})$. Equipped with this estimate, the
trivial bound $F(\al)\le (2P+1)^s$ suffices to conclude that
$$ \Ups_0 \ll P^{s +(55\dlt k -k)(s-1)} \ll P^{s-2k}. $$
On collecting together the various estimates, the lemma then follows from
\rf{A}.

We address the set $M_\al(C)$ by first investigating the consequences
of the defining conditions $|f(\lm_j\al)|\ge P^{1-11\dlt}$ individually.
One  applies Weyl's inequality (\cite{hlm}, Lemma 2.4) in reverse,
followed by a joint application of Theorems 4.1 and 2.8 of \cite{hlm}.
This shows that there exists some real number $K\ge 1$ (depending only on $k$),
 and integers
$a_j,q_j$ with $1\le q_j\le KP^{11k\dlt}$ and
$$ |\lm_j\al -a_j/q_j | \le q_j^{-1} KP^{11k\dlt -k}. $$
Define real numbers $\theta_j$ by
$$  \lm_j\al = \f{a_j}{q_j} + \theta_j q_j^{-1}P^{11k\dlt -k}. $$
Then $|\theta_j|\le K$. Also, when $\bflm\in M_\al(C)$, then $|\lm_j\al|\ge
C^{-1}P^{1/3-k}$, and $11k\dlt < \f13$ holds for all $k\ge 2$. Hence,
when $P$ is sufficiently large in terms of $k$ and $C$, we must have 
$a_j\neq 0$. It follows that
$$ \f{\lm_2}{\lm_1} = \f{\lm_2\al}{\lm_1\al} = \f{q_1 a_2}{q_2 a_1}
\Big(1+ \f{\theta_2}{a_2} P^{11k\dlt -k}\Big)
\Big(1+ \f{\theta_1}{a_1} P^{11k\dlt -k}\Big)^{-1}. $$
In particular, when $P$ is large, one infers that
$$ \Big| \f{\lm_2}{\lm_1} - \f{q_1 a_2}{q_2 a_1}\Big| \le 4KC^2 P^{11k\dlt -k}. $$
Now write $q_1a_2/(q_2 a_1) = a/q$ in lowest terms. Then $q|q_2 a_1$, and
$$|a_1| \le |\lm_1 \al q_1| + KP^{11k\dlt -k} \le 2CKP^{11k\dlt}. $$
Consequently,
\be{q} q\le 2CK^2 P^{22k\dlt}. \ee
We have now shown that whenever $\al$ is in the indicated range, and $\bflm
\in M_\al(C)$, then $\lm_2/\lm_1 \in \fk K$ where $\fk K$ denotes the
union of all intervals $\{\lm: \, |\lm - a/q| \le 4KC^2 P^{11k\dlt-k}\}$ with
$a\in\ZZ$, $(a,q)=1$ and $q\in\NN$ satisfying  \rf{q}. Also, 
$ |\lm_2/\lm_1|\le C^2$, and we therefore write ${\fk K}(C)= {\fk K}\cap
[-C^2, C^2]$. By symmetry, we now have $\lm_j/\lm_1\in {\fk K}(C)$ for
all $2\le j\le s$.
Since the measure of  ${\fk K}(C)$ is $O(P^{55k\dlt-k})$, the transformation
$\mu_j=\lm_j/\lm_1$ yields
$$ \int_{M_\al(C)} \,d\bflm \ll \int_{1/C}^C \lm_1^{s-1}\,d\lm_1 
\Big( \int_{{\fk K}(C)} \,d\mu\Big)^{s-1} \ll P^{(55k\dlt-k)(s-1)}, $$
as required. This completes the proof of Lemma 6.

\section{The principal proposition}
We are ready to derive a variant of the theorem in which
$N_{\bflms} (P,\tau)$ is replaced by $I_{\bflms}(P,\eta)$. Let $P\ge 2$ and
$0<\eta\le 2$. Then define ${\cal K}(P,\eta)$ to be the set of  all
$\bflm\in\Lm_k^{(s)}(C)$ where
\be{71} |I_{\bflms}(P,\eta) - \eta J_{k,s}(\bflm) P^{s-k} | >2P^{s-k-3\dlt}. \ee
Moreover, let ${\cal K}(P)$ be the union of all ${\cal K}(P,\eta)$ with
$0<\eta\le 2$. Note that the left hand side of \rf{71} defines a measurable
function of $\bflm$ so that ${\cal K}(P,\eta)$ and ${\cal K}(P)$ are
measurable sets.  The first step in the argument is to show that
\be{72} \int_{{\cal K}(P)} \,d\bflm \ll P^{-3\dlt}. \ee

It suffices to prove \rf{72} for 
$P\ge P_0(C)$ where $P_0$ is a suitably large number depending only on $C$.
For technical convenience, we consider a slight modification of 
${\cal K}(P,\eta)$. Let $P\ge P_0$ and $0<\eta\le 2$, and let 
${\cal H}(P,\eta)$  denote the set of all $\bflm \in \Lm_k^{(s)}(C)$ where
\be{72a} |I_{\bflms}(P,\eta) - \eta J_{k,s}(\bflm) P^{s-k} | >P^{s-k-3\dlt}. \ee
To obtain an estimate for the measure of ${\cal H}(P,\eta)$, we apply Lemma 5 
and recall \rf{12a}. Then \rf{72a} yields
$$|{\cal I}({\fk c})| > P^{s-k-3\dlt} - O(P^{s-k-1/3}) \ge \f12 P^{s-k-3\dlt},  
$$
because $P$ is large. Consequently, by \rf{12} and Lemma 6, 
\be{73} \int_{{\cal H}(P,\eta)} \,d\bflm \le 2P^{3\dlt+k-s} \int_{\Lm_k^{(s)}(C)}
|{\cal I}({\fk c})|\, d\bflm \ll P^{-7\dlt}.\ee

Now let $L=[P^{4\dlt}]$. For $1\le l\le 2L$, put $\eta_l = l/L$. First 
suppose that $\bflm \in {\cal K}(P,\eta)$, for some $\eta\in[1/L, 2]$. Then,
there is some $l$ with $1\le l\le 2L$ and $\eta_l\le \eta \le \eta_{l+1}$.  
By \rf{71}, one of the inequalities
\begin{eqnarray}
I_{\bflms}(P,\eta) - \eta J_{k,s}(\bflm) P^{s-k} & >& 2P^{s-k-3\dlt},\label{74}
 \\
 \eta J_{k,s}(\bflm) P^{s-k}- I_{\bflms}(P,\eta) & >& 2P^{s-k-3\dlt} \label{75}
\end{eqnarray}
must hold. The function $I_{\bflms}(P,\eta)$ is increasing, as a 
function of $\eta$. Hence, if \rf{75} holds,
then
\begin{eqnarray*}
I_{\bflms}(P,\eta_l) &\le& I_{\bflms}(P,\eta) < \eta J_{k,s}(\bflm)P^{s-k}
- 2 P^{s-k-3\dlt} \\
& \le & \eta_l  J_{k,s}(\bflm)P^{s-k}- 2 P^{s-k-3\dlt} + O(L^{-1}P^{s-k}),
\end{eqnarray*}
with an implicit constant depending only on $C$. Hence, for sufficiently 
large $P$, it follows that \rf{72a} holds with $\eta=\eta_l$, and hence that
$\bflm\in{\cal H}(P,\eta_l)$. Similarly, if \rf{74} holds, one uses
that $I_{\bflms}(P,\eta)\le I_{\bflms}(P,\eta_{l+1})$, and finds that  
$\bflm\in{\cal H}(P,\eta_{l+1})$. This shows that the union of 
${\cal K}(P,\eta)$ with $1/L\le \eta\le 2$ is contained in the union of the
${\cal H}(P,\eta_l)$ with $1\le l\le 2L$, and by \rf{73}, the measure 
of this set does not exceed $O(LP^{-7\dlt}) = O(P^{-3\dlt})$, as required.

Now suppose that $\bflm \in {\cal K}(P,\eta)$, for some $\eta$ with $0<\eta
\le 1/L$. Then $$\eta J_{k,s}(\bflm)P^{s-k}=O(P^{s-k-4\dlt}),$$ with the 
implicit constant depending only on $C$. Thus, when $P$ is 
sufficiently large, \rf{71} 
implies that 
$$ I_{\bflms} (P,\eta_1) \ge I_{\bflms} (P,\eta)> \f32P^{s-k-3\dlt}. $$
One may subtract $\eta_1 J_{k,s}(\bflm)P^{s-k}$ here. By a now familiar 
reasoning, one then finds that \rf{72a} holds with $\eta=\eta_1$. Hence,
the union of ${\cal K}(P,\eta)$ with $0<\eta\le 1/L$ is contained in 
${\cal H}(P,\eta_1)$, and so, by \rf{73}, has measure not exceeding 
$O(P^{-7\dlt})$. Combining this with the previous discussion, one
confirms that \rf{72} indeed holds.

\begin{lem} Let $s>2k\ge 4$ and $\dlt=8^{-2k}$. Then, for almost all 
$\bflm\in\Lm_k^{(s)}$, there exists a number $P_0=P_0(\bflm,k)$ such that the
inequality
\be{77}
 |I_{\bflms}(P,\eta) - \eta J_{k,s}(\bflm) P^{s-k} | \le P^{s-k-5\dlt/2} 
\ee
holds for all $P\ge P_0$ and $0<\eta\le 2$.
\end{lem}

{\em Proof}. Let $r=3/(8\dlt)$. Then $r\in\NN$, and the series 
$\sum_j j^{-3r\dlt}$ converges. Hence, by \rf{72}, for any $\eps>0$,
any $C\ge 2$, there is a number $J=J(\eps,C)$ such that the set
$$ {\cal K}_\eps = \bigcup_{j\ge J} {\cal K}(j^r) $$
has measure not exceeding $\eps$.  
For $\bflm\in\Lm_k^{(s)}\setminus{\cal K}_\eps$, one then has
\be{76}
 |I_{\bflms}(P,\eta) - \eta J_{k,s}(\bflm) P^{s-k} | \le 2P^{s-k-3\dlt}
\ee
for all $P=j^r$ with $j\ge J$, uniformly in $0<\eta \le 2$.

Now suppose that $\bflm\in\Lm_k^{(s)}(C)$ is an $s$-tuplet where a number
$P_0$ with the properties described in Lemma 7 does {\em not} exist.
We have just seen that this is possible only when $\bflm\in {\cal K}_\eps$
holds for {\em all} $\eps>0$. In particular, $\bflm$ is in the intersection of 
 ${\cal K}_{1/m}$ with $m\in\NN$, and the latter is a null set. Hence, the
set of all $\bflm\in\Lm_k^{(s)}(C)$ where a $P_0$ as desired does not exist, also has Lebesgue measure zero. Now consider the union of these sets, with $C
\ge 2$ running over natural numbers, to complete the proof of Lemma 7.

\section{The final sandwich}

In this section we remove the weights from the counting function $I_{\bflms}(P,\eta)$ by a standard sandwich technique.
Let $0<\tau\le 1$, and put $\Delta=P^{-5\dlt/4}$. Consider the functions $W^+$,
$W^-$, defined by
\be{81} W^\pm = \pm \Delta^{-1} \big((\tau\pm \Delta) W_{\tau\pm \Delta} (\al) -\tau W_\tau(\al)\big). \ee
The $W^\pm$ are continuous functions with $0\le W^{\pm}(\al)\le 1$
for all $\al\in\RR$, and one readily checks that
$$\begin{array}{llll}
W^+ (\al) = 1 & \mbox{for $|\al|\le\tau$}, & W^+(\al)=0 & \mbox{for $|\al|\ge \tau+\Delta$}, \\
W^- (\al) = 1 & \mbox{for $|\al|\le\tau-\Delta$,\;\;} & W^+(\al)=0 & \mbox{for $|\al|\ge \tau$.}
\end{array} 
$$
Now define the counting functions
$$ I^\pm(P,\tau)= \sum_{|\bfx|\le P} W^\pm(\lm_1 x_1^k + \ldots + \lm_s x_s^k), $$
and note that one has the sandwich inequalities
\be{82} I^-_{\bflms}(P,\tau) \le N_{\bflms}(P,\tau) \le  
I^+_{\bflms}(P,\tau). \ee
By \rf{81}, we have
$$ I^\pm(P,\tau)= \pm \Delta^{-1} \big((\tau\pm \Delta) I_{\bflms}(P,\tau\pm \Delta) -\tau I_{\bflms}(P,\tau)\big). $$
Now suppose that $\bflm$ is not in the exceptional null set excluded from the claim in Lemma 7. Then, for large $P$, one has
\begin{eqnarray*}
 I^\pm(P,\tau) &=&  \pm \Delta^{-1} \big( ((\tau\pm\Delta)^2-\tau^2)J_{k,s}(\bflm)P^{s-k} + O((\tau+\Delta)P^{s-k-5\dlt/2})\big) \\
& = & 2\tau J_{k,s}(\bflm)P^{s-k}+ O(\Delta P^{s-k} + \tau \Delta^{-1}P^{s-k-5\dlt/2}).
\end{eqnarray*}
The error term here is $O(P^{s-k-5\dlt/4})$, uniformly in $\tau$. Hence, by 
\rf{82}, 
$$ |N_{\bflms}(P,\tau) - 2\tau J_{k,s}(\bflm)P^{s-k}| < P^{s-k-\dlt}. $$
This proves that if $\bflm$ is in the exceptional set in the theorem, then
$\bflm$ is also in the exceptional set in Lemma 7. This completes the proof of the theorem.

\end{document}